\documentclass{article}
\usepackage[pdftex]{hyperref}
\hypersetup{
  colorlinks   = true, 
  urlcolor     = gray, 
  linkcolor    = red, 
  citecolor   =  blue
}
\usepackage{amsmath,amsfonts,amsthm,amssymb,graphicx,xcolor}
\usepackage[all,2cell,ps]{xy}

\theoremstyle{plain}
\newtheorem{thm}{Theorem}[section]
\newtheorem{lem}[thm]{Lemma}
\newtheorem{prop}[thm]{Proposition}
\newtheorem{cor}[thm]{Corollary}

\theoremstyle{definition}
\newtheorem{defn}[thm]{Definition}

\newtheorem{rem}[thm]{Remark}
\newtheorem{rems}[thm]{Remarks}

\newcommand{\Z}{\mathbb Z}

\newcommand{\C}{\mathbb C}
\newcommand{\pp}{\mathbb{P}}

\DeclareMathOperator{\Hom}{Hom}
\DeclareMathOperator{\Ext}{Ext}

\newcommand{\ssm}{\smallsetminus}

\newenvironment{pf}{\begin{proof}}{\end{proof}}

\title{Exceptional collections and the bicanonical map of Keum's fake projective planes}
\author{\small{Gennaro Di Brino}\footnote{Partially supported by the National Research Fund, Luxembourg, and the Marie Curie Actions Program of the European Commission (FP7-COFUND).} \\ \scriptsize{Max Planck Institute for Mathematics\footnote{Most recent institutional affiliation.}}\\ \footnotesize{\textsf{gennaro.dibrino@gmail.com}} \and \small{Luca F. Di Cerbo}\footnote{Partially supported by a grant of the Max Planck Society: ``Complex Hyperbolic Geometry and Toroidal Compactifications'', and by a grant associated to the position at ICTP.} \\ \scriptsize{International Centre for Theoretical Physics} \\ \footnotesize{\textsf{ldicerbo@ictp.it}}}
\date{}

\begin{document}

\maketitle

\begin{abstract}
We apply the recent results of Galkin \emph{et al}. \cite{Katzarkov} to study some geometrical features of Keum's fake projective planes. Among other things, we show that the bicanonical map of Keum's fake projective planes is always an embedding. Moreover, we construct a nonstandard exceptional collection on the unique fake projective plane $X$ with $H_{1}(X; \Z)=(\Z/2\Z)^{4}$.
\end{abstract}

\section{Introduction}

This paper concerns the geometry of Keum's fake projective planes.  Recall that Keum's fake projective planes are those planes with automorphism group $G_{21}$, where by $G_{21}$ we denote the non abelian group of order $21$. We refer to Section \ref{FPP} for more details. The first result presented in this paper studies the bicanonical map of such planes. More precisely, we have the following.

\begin{thm}\label{first}
Let $X$ be a fake projective plane with $Aut(X)=G_{21}$. The $2$-canonical map 
\[
\varphi_{|2K_{X}|}: X\rightarrow \pp^{9}
\]
is an embedding.
\end{thm}
 
Along the way, we also show that the bicanonical map of \emph{any} fake projective plane is a birational morphism. This result was originally proved by Mendes-Lopes and Pardini in \cite{Pardini}. The proof presented here is of a different flavor and the final result somewhat stronger.

Furthermore, building on the work of Galkin \emph{et al}. \cite{Katzarkov}, we study exceptional collections on Keum's fake projective planes. In \emph{loc. cit.}, it is shown that any of Keum's planes supports a  ``standard" exceptional collection. Here we show that one of Keum's planes admits a \emph{nonstandard} exceptional collection, i.e., a collection which does not coincide with the one already constructed in \cite{Katzarkov}. More precisely, we have the following.

\begin{thm}\label{second}
Let $X$ be the fake projective plane with $Aut(X)=G_{21}$ and with $H_{1}(X; \Z)=(\Z/2\Z)^{4}$. We have a nonstandard exceptional collection.
\end{thm}

For a more precise statement we refer to Theorem \ref{eccezionale} in Section \ref{exotic}. The proof of Theorem \ref{second} relies on a vanishing result for certain $(\Z/7\Z)$-equivariant line bundles on Keum's planes, see Proposition \ref{vanishing} in Section \ref{exotic}. This vanishing result complements the vanishing theorem for $G_{21}$-equivariant line bundles presented in \cite{Katzarkov}. The techniques used in the proof of this result are distinct from the ones used in \emph{loc. cit.}. Finally, we end our treatment with a translation of our geometrical results into the language of derived categories of coherent sheaves and semi-orthogonal decompositions. For more details see Corollary \ref{orthogonal} in Section \ref{exotic}. \\

Most of the results presented here crucially rely on the detailed analysis presented by Keum in \cite{Keum}, of the action of the automorphism group on fake projective planes.\\

\noindent\textbf{Acknowledgements}. The authors thank Prof. Fabrizio Catanese for sharing his knowledge with them, for suggesting Proposition \ref{catanese} and for improving an earlier version of Theorem \ref{birat}. The first named author would like to thank Prof. Mikhail Kapranov for helpful conversations and Prof. Yuri Manin for his kind interest in this work. The authors thank the Max Planck Institute for Mathematics and the Mathematics Research Unit of the University of Luxembourg for their hospitality at the beginning of this project. Finally, the authors thank the Insitut des Hautes \'Etudes Scientifiques, the Max Planck Institute for Mathematics and the International Centre for Theoretical Physics for the excellent working environments.

\section{Preliminaries}\label{FPP}

In this section, we recall the general properties of the so-called fake projective planes. We briefly describe a particular class of highly symmetrical fake projective planes first studied by Keum. The results and definitions collected in this section are used throughout the rest of this work.

\begin{defn}
A fake projective plane is a surface of general type $X$ with $c_{2}=3$ and $p_{g}=H^{0}(X; K_{X})=0$.
\end{defn}

Interestingly, any such surface has the same Hodge diamond as $\pp^{2}$. A priori it is not at all clear that any such surface exists.  It turns out that it is quite hard to construct such objects, and the first fake projective plane was obtained by Mumford in \cite{Mumford} via a highly nontrivial construction. Because of Yau's solution to Calabi's conjecture \cite{Yau}, it is well known that any fake projective plane is necessarily a compact complex hyperbolic $2$-manifold. In other words, any fake projective plane $X$ is given as the quotient of the complex hyperbolic $2$-space $\mathcal{H}^{2}_{\C}$ by a torsion free co-compact arithmetic lattice $\Gamma\in \text{PU}(2, 1)$. The arithmeticity of the fundamental group of a fake projective plane is a consequence of the works of Klingler \cite{Klingler} and Yeung \cite{Yeung}.

The fundamental groups of fake projective planes have been classified by Prasad-Yeung \cite{Prasad} and Cartwright-Steger \cite{Steger}. More precisely, there are exactly 50 lattices in $\text{PU}(2, 1)$ which are fundamental groups of fake projective planes. This computer assisted classification not only explicitly provides the possible fundamental groups, but also computes the abelianization and the automorphism group for any such $\Gamma\in \text{PU}(2, 1)$. It turns out that (up to complex conjugation) there are exactly three fake projective planes with automorphism group isomorphic to $G_{21}$, where $G_{21}$ denotes the unique non abelian group of order $21$. All of the other admissible automorphism groups have smaller cardinality and in many cases are just reduced to the identity. Fake projective planes with automorphism group $G_{21}$ are called Keum's planes as the first such example was described by Keum in \cite{Keum7}.

Let us observe that, given a fake projective plane $X$ with $\pi_{1}(X)=\Gamma$, we have that
\[
H_{1}(X; \Z)=\Gamma/[\Gamma, \Gamma]
\]
is a never vanishing torsion group. Recall that since $H^{1}(X; \mathcal{O})=0$, we have that the group $H_{1}(X; \Z)$ must always be torsion. For the explicit list of possible first homology groups, we refer once again to \cite{Steger}. 

Next, let us collect a few useful facts regarding the Picard group of a fake projective plane. The vanishing $H^{1}(X; \mathcal{O})=0$ implies that the Picard group $Pic(X)$ is always isomorphic the the cohomology group $H^{2}(X; \Z)$. Thus, by the universal coefficient theorem we compute that $\text{Tor}(H^{2}(X; \Z))=H_{1}(X; \Z)\neq 0$. In other words, on any fake projective plane there are torsion line bundles. Since this fact will be needed in the rest of this paper, we summarize this discussion into a lemma.

\begin{lem}\label{Tor}
For a fake projective plane X, we have $Pic(X)\simeq H^{2}(X; \Z)$. Moreover, the torsion part of $H^{2}(X; \Z)$ is equal to $H_{1}(X; \Z)$ which is never vanishing. 
\end{lem}

Now, by Poincar\'e duality $Pic(X)/\text{Tor}(H^{2}(X;\Z))$ is a one dimensional unimodular lattice which is generated by an ample line bundle $L$ such that $c_{1}(L)^{2}=1$. Thus, it is worth making the following definition.

\begin{defn}
For a fake projective plane X, we denote by $L_{1}$ any ample generator of the torsion free part of $Pic(X)$. 
\end{defn}

Of course, the line bundle $L_{1}$ is not uniquely determined. In fact, the different choices are parametrized by the torsion line bundles on $X$, i.e., they are in one-to-one correspondence with $H_{1}(X; \Z)$. 

Let us now denote by $``\equiv"$, resp. ``$\cong$'', the relation of numerical equivalence, resp. linear equivalence, for line bundles and divisors. Since a fake projective plane has Picard number equal to one, given any ample line bundle $L$ with self-intersection $k^{2}$ we have that
\[
L\equiv kL_{1},
\]
where $L_{1}$ is any ample generator of the torsion free part of $Pic(X)$. For our purposes, it is then natural to make the following definition.

\begin{defn}
For a fake projective plane X, we denote by $L_{k}$ any ample line bundle $L$ such that $c^{2}_{1}(L)=k^{2}$. 
\end{defn}

We end this section by briefly describing the notion of an exceptional collection. 

\begin{defn}
Let $\mathcal{D}^b(X)$ be the bounded derived category of coherent sheaves on a smooth projective variety $X$. A collection of objects $\left(E_0,\ldots,E_{r-1}\right)$ in $\mathcal{D}^b(X)$ is said to be \emph{exceptional} if $\Ext^k(E_j,E_i)$ vanishes whenever $j>i$ and for all $k\in\Z$, whereas $\Hom(E_i,E_i)$ has rank one for $i=0,\ldots, r-1$. We will denote the $\Ext$'s by $\Hom(E_j,E_i[k])$ when no confusion can arise.
\end{defn}

 Given an exceptional collection, one may ask whether its triangulated envelope, that is, the smallest triangulated subcategory of $\mathcal{D}^b(X)$ containing such a collection, coincides with $\mathcal{D}^b(X)$. The collection is called \textit{full} if this is the case. When a given collection is not full (see Remark \ref{R:exotica}), one can consider its ``orthogonal complement'' $\mathcal{A}$ which then, together with the exceptional collection, generates $\mathcal{D}^b(X)$. From a more categorical viewpoint, it also makes sense to ask whether $\mathcal{D}^b(X)$ admits a so-called \emph{phantom} in its Hochschild homology. That is, a nontrivial subcategory $\mathcal{A}$ of $\mathcal{D}^b(X)$ such that $HH_{\bullet}(\mathcal{A})$ vanishes. We address this question in our case in Remark \ref{R:phantom}.\\

For further details and a good review of the literature on derived categories, we refer the reader to \cite{Bondal},  \cite{Katzarkov}, \cite{Kuz} and the references therein. On the other hand, for the basic complex surface theory and complex hyperbolic geometry, the reader is referred to \cite{Bar}, \cite{Bea} and \cite{Gold Book}.\\

\section{Bicanonical map of Keum's fake projective plane}

In this section, we study the bicanonical map of fake projective planes. More precisely, we focus on the bicanonical map of Keum's fake projective planes. Let us begin with a series of elementary lemmas regarding the cohomology of line bundles over fake projective planes. Recall that, on any fake projective plane $X$, we denote by $L_{k}$ any line bundle numerically equivalent to $kL_{1}$, where $L_{1}$ is any generator of $Pic(X)$ and $k\geq 1$ any integer. 

\begin{lem}\label{L4}
For any fake projective plane X, we have $h^{0}(X; L_{4})=3$.
\end{lem}
\begin{pf}
First, let us observe that
\[
\chi(X; L_{4})=h^{0}(X; L_{4})-h^{1}(X; L_{4})+h^{2}(X; L_{4})=1+\frac{1}{2}4L_{1}\cdot(4L_{1}-K_{X})=3.
\]
Recall that $K_{X}\equiv L_{3}$, so that $L_{4}\equiv K_{X}+L_{1}$. By Kodaira vanishing, we have $h^{1}(X; L_{4})=h^{2}(X; L_{4})=0$. Thus, for any for any $L_{4}$ on $X$, we have the identity $h^{0}(X; L_{4})=3$. 
\end{pf}

We now have the following estimate for the cohomology of a $L_{2}$ on $X$.

\begin{lem}\label{L2}
For any fake projective plane $X$ and for any $L_{2}$, we have $h^{0}(X; L_{2})\leq 2$.
\end{lem}
\begin{pf}
By contradiction, let us assume that $h^{0}(L_{2})>2$. Because of Lemma 15.6.2 in \cite{Kol}, the multiplication map
\[
H^{0}(X; L_{2})\times H^{0}(X; L_{2})\rightarrow H^{0}(X; L_{4})
\]
is such that $h^{0}(X; L_{2}+L_{2})\geq h^{0}(X; L_{2})+h^{0}(X; L_{2})-1\geq 5$. This fact contradicts Lemma \ref{L4}.
\end{pf}

The same reasoning applies to the study of the cohomology of any $L_{1}$ as well. More precisely, we have the following.

\begin{lem}\label{L1}
Let $X$ be a fake projective plane. For any $L_{1}$, if $h^{0}(L_{1})\neq 0$, we then have $h^{0}(X; L_{1})=1$.
\end{lem}
\begin{pf}
First, if we assume $h^{0}(L_{1})\geq 2$, the multiplication map
\[
H^{0}(X; L_{1})\times H^{0}(X; L_{1})\rightarrow H^{0}(X; L_{2})
\]
gives that $h^{0}(X; L_{2})\geq 3$. Then, we can finish by using Lemma \ref{L2}. 
\end{pf}

We can then formulate the following useful lemma, which gives the finiteness of curves numerically equivalent to $L_{1}$ inside a fake projective plane.

\begin{lem}\label{Curves}
Let $X$ be a fake projective plane. There are at most finitely many curves $C_{1}, ..., C_{k}$ in $X$ numerically equivalent to $L_{1}$.
\end{lem}
\begin{pf}
Let us observe that if there exists an effective $L_{1}$ on $X$, this must necessarily correspond to a unique curve, say $C_{1}$, in $X$. This fact follows from Lemma \ref{L1}. Now, the number of different $L_{1}$'s on $X$ is parametrized by $H_{1}(X; \Z)$ which is a finitely generated torsion group, see Lemma \ref{Tor}. Thus, there are at most finitely many $C_{i}$'s numerically equivalent to $L_{1}$. The proof is complete.
\end{pf}

\begin{rem}\label{smooth}
It seems currently unknown whether or not there exists an effective $L_{1}$ on any of the fake projective planes. Nevertheless, it can be proved that if a curve as in Lemma \ref{Curves} exists, then it is necessarily smooth of genus three. See Corollary 2.4 in \cite{ultimo}.
\end{rem}

The finiteness result given in Lemma \ref{Curves} and Remark \ref{smooth} have nice consequences for the birational geometry of fake projective planes. More precisely, we have the following.

\begin{thm}\label{birat}
Let $X$ be a fake projective plane. The $2$-canonical map 
\[
\varphi_{|2K_{X}|}: X\rightarrow \pp^{9}
\]
is a birational morphism, and an isomorphism with its image outside a finite set of points in $X$.
\end{thm}
\begin{pf}
First, let us observe that 
\[
h^{0}(X; 2K_{X})=1+K^{2}_{X}=10
\]
so that the bicanonical map goes into $\pp^{9}$ as claimed. Next, let us observe that in a fake projective plane $X$, there are no curves $C$ such that $C^{2}=0, -1$. Since $K_{X}$ is ample with $K^{2}_{X}=9$, the first part of Reider's theorem (see page 176 in \cite{Bar}) tells us that $\varphi_{|2K_{X}|}$ is indeed a morphism. Because of Lemma \ref{Curves}, there are at most finitely many curves, say $C_{1}, ..., C_{k}$, numerically equivalent to $L_{1}$. The second part of Reider's theorem ensures that if two points $P$ and $Q$, possibly infinitely near, are not separated by $\varphi_{|2K_{X}|}$, then they are both contained in at least one of the $C_{i}$'s.
Thus, the map $\varphi_{|2K_{X}|}$ is a biholomorphism on $X\ssm C_{1}\cup ...\cup C_{k}$. The proof of the birationality of the morphism $\varphi_{|2K_{X}|}$ is then complete.

For the remaining part of the statement we argue as follows. Let $C$ be any curve numerically equivalent to $L_{1}$.  We can then write
\[
K_{X}\cong 3C+T
\]
where $T\in Pic(X)$ is a \emph{torsion} line bundle. Next, let us observe that because of Kodaira vanishing, the restriction map
\[
H^{0}(X; 2K_{X})\rightarrow H^{0}(C; (2K_{X})_{|_{C}})
\] 
is surjective. Also, by adjunction, we have 
\[
(2K_{X})_{|_{C}}\cong K_{C}\otimes \mathcal{O}_{C}(2C+T).
\]
Thus, given two points $P, Q\in C$, possibly infinitely near, let $\mathcal{I}_{P, Q}$ be their ideal sheaf in $\mathcal{O}_C$. Consider then the following short exact sequence of sheaves on $C$
\[
0\rightarrow \mathcal{I}_{P, Q}\otimes K_{C}\otimes \mathcal{O}_{C}(2C+T)\rightarrow  K_{C}\otimes \mathcal{O}_{C}(2C+T)\rightarrow \C^{2}\rightarrow 0.
\]
By \cite{Cat96} and \cite{Cat99} we have that the points $P, Q$ are not separated by $ K_{C}\otimes \mathcal{O}_{C}(2C+T)$ if and only if 
\[
H^{1}(C; \mathcal{I}_{P, Q}\otimes K_{C}\otimes \mathcal{O}_{C}(2C+T))=\C.
\]
By duality, the points $P, Q$ are not separated if and only if 
\[
\mathcal{I}_{P, Q}\otimes \mathcal{O}_{C}(2C+T)\simeq\mathcal{O}_{C}.
\]
Thus, the sheaf $\mathcal{I}_{P, Q}$ is invertible and $P+Q$ is the \emph{unique} divisor of a section in $H^{0}(C; \mathcal{O}_{C}(2C+T))\simeq\C$. Concluding, by Lemma \ref{Curves} there are at most finitely many curves $C$ numerically equivalent to $L_{1}$, and then the birational morphism $\varphi_{|2 K_{X}|}$ is necessarily an isomorphism with its image outside a finite set of points in $X$.
\end{pf}

\begin{rems} {\bf (1)} The reader will notice that the above argument also shows that the bicanonical map is a finite morphism.\\
{\bf (2)} The birationality of $\varphi_{|2K_{X}|}$ for fake projective planes was originally proved by M. Mendes-Lopes and R. Pardini with a different argument, see Theorem 1.1 in \cite{Pardini}.
\end{rems}

Next, let us show that Theorem \ref{birat} can be improved for the so-called Keum's fake projective planes. More precisely, using a vanishing theorem first proved in \cite{Katzarkov}, it can be proved that the bicanonical map of Keum's fake projective planes is actually an embedding. Following \emph{loc. cit.}, Keum's fake projective planes are those planes with automorphism group $G_{21}$.  According to the list given by Cartwright and Steger \cite{Steger}, there are exactly 3 distinct such planes (up to complex conjugation).

\begin{thm}\label{Keum}
Let $X$ be a fake projective plane with $Aut(X)=G_{21}$. The $2$-canonical map 
\[
\varphi_{|2K_{X}|}: X\rightarrow \pp^{9}
\]
is an embedding.
\end{thm}

\begin{pf}
Let $X$ be a fake projective plane with automorphism group $G_{21}$. By Lemma 2.2 in \cite{Katzarkov}, there exists a unique $G_{21}$-equivariant line bundle, say $\mathcal{O}_{X}(1)$, such that $K_{X}\cong 3\mathcal{O}_{X}(1)$. Recall that  $``\cong"$ denotes linear equivalence. The main technical result in \emph{loc. cit.} (see Theorem 1.3 therein) is the following remarkable vanishing:
\[
H^{0}(X; 2\mathcal{O}_{X}(1))=0.
\]
This vanishing has some interesting consequences. Recall that Cartwright and Steger explicitly computed $H_{1}(X;
\Z)$ for any fake projective plane $X$.  In our case, if $X$ is a fake projective plane with $Aut(X)=G_{21}$, we have
\[
H_{1}(X; \Z)=(\Z/2\Z)^{3},\quad (\Z/2\Z)^{4},\quad \text{or} \quad (\Z/2\Z)^{6}.
\]
Thus, for any torsion line bundle $T\in Pic(X)=H^{2}(X; \Z)$, we have $2T\cong \mathcal{O}_{X}$. Given any $L_{1}$ on $X$, we can write $L_{1}\cong \mathcal{O}_{X}(1)+T$ for some torsion line bundle $T\in Pic(X)$. This implies that given any $L_{1}$ on $X$, we necessarily have $2L_{1}\cong 2\mathcal{O}_{X}(1)$. We then conclude that
\[
H^{0}(X; L_{1})=0
\]
for any $L_{1}$ on $X$. In other words, if $X$ is a fake projective plane with $Aut(X)=G_{21}$, there are no curves $C_{i}$ numerically equivalent to $L_{1}$. Since on $X$ there are clearly no curves $C$ with $C^{2}=0, -1$, or  $-2$, the second part of Reider's theorem tells us that, if $\varphi_{|2K_{X}|}$ is not an embedding, there has to be at least one curve numerically equivalent to $L_{1}$. Since we know that no such curves can exist on $X$, the argument is complete.
\end{pf}

\begin{rem}\label{rimarco}
It is certainly tempting to conjecture that there are no effective $L_{1}$ on any of the fake projective planes. If this conjecture holds true, the argument given in Theorem \ref{Keum} would imply that $\varphi_{|2K_{X}|}$ is always an embedding.
\end{rem}

It seems of interest to conclude this section with a result precisely addressing when the bicanonical map of a fake projective plane is an \emph{embedding}. As pointed out in Remark \ref{rimarco}, the non-existence of an effective $L_{1}$ certainly would suffice. Nevertheless, a somewhat weaker requirement gives a necessary and sufficient condition for $\varphi_{|2K_{X}|}$ to be an embedding. More precisely, we have the following.

\begin{prop}\label{catanese}
Let $X$ be a fake projective plane. The $2$-canonical map 
\[
\varphi_{|2K_{X}|}: X\rightarrow \pp^{9}
\]
is an embedding if and only if there is no effective divisor $C$ with $C^{2}=1$ and $h^{0}(C; \mathcal{O}_{C}(2C))=h^{1}(C; \mathcal{O}_{C}(2C))=1$.
\end{prop}

\begin{pf}
In the proof of Theorem \ref{birat}, we have seen that $\varphi_{|2K_{X}|}$ is an embedding if and only if there are no effective divisors $C$ with $C^{2}=1$ and such that
\[
H^{0}(C; \mathcal{O}_{C}(2C+T))\simeq \C
\] 
where $T\in Pic(X)$ is the torsion line bundle such that $K_{X}\cong 3C+T$. By Riemann-Roch, we have
\[
h^{0}(C;  \mathcal{O}_{C}(2C+T))=h^{1}(C;  \mathcal{O}_{C}(2C+T))
\]
and by duality
\[
h^{1}(C;  \mathcal{O}_{C}(2C+T))=h^{0}(C; K_{C}\otimes\mathcal{O}_{C}(-2C-T))=h^{0}(C; \mathcal{O}_{X}(2C))
\]
so that
\[
h^{0}(C; \mathcal{O}_{X}(2C))=h^{1}(C; \mathcal{O}_{C}(2C))=1.
\]
The proof is then complete.
\end{pf}

\section{A nonstandard exceptional collection}\label{exotic}

In this section, we will again assume that the surface $X$ is a fake projective plane, so that $Pic(X)$ is of rank one and $K_X^2=9$. The statement below follows from the definition of exceptional collection combined with some complex surface theory, and it has already appeared in the recent literature (\cite{Katzarkov}, \cite{Keumvan}, \cite{LY}). We give the details for the reader's convenience.

\begin{prop}\label{P:cubicroot}
Let $X$ be as above and let $L'$ be an ample generator of $Pic(X)$. Then the sequence $\left(\mathcal{O}_{X}, -L', -2L'\right)$ is exceptional if and only if the dimensions $h^{0}(X; 2L')$, $h^{2}(X; L')$ and $h^{2}(X; 2L')$ vanish.
\end{prop}

\begin{pf}
We start by noticing that for the bundle $L'$ in the statement, being an ample generator of $Pic(X)$ is equivalent to being a numerical cubic root of the canonical bundle, that is, $3L'\equiv K_X$. This said, we check that the sequence in the statement satisfies the definition of an exceptional collection. First, notice that $h^{0}(X; 2L')=0$ implies $h^{0}(X; L')=0$ for elementary reasons. Next, let $M$ denote either $L'$ or $2L'$, and let $D$ be the divisor associated with $M$. From the Hirzebruch-Riemann-Roch formula one gets
\begin{equation}\label{e:fakeHRR}
\chi(M)=\frac{1}{2}D\cdot(D-K_X)+1-q+p_{g}
\end{equation}
where $q(X)$ and $p_{g}(X)$ are respectively the irregularity and geometric genus of the surface $X$. Since $X$ is a fake projective plane, we have $p_{g}(X)=q(X)=0$ and, by assumption, $h^2(X; M)=0$. Expanding the Euler characteristic on the left, we have that $-h^1(M)$ is equal to the right hand side of equation \eqref{e:fakeHRR}. But the latter is in turn zero by direct computation, for $M$ equal to either $L'$ or $2L'$.\\

Let now $\left(E_0, E_1, E_2\right)=\left(\mathcal{O}_{X}, -L', -2L'\right)$, and let $k$ denote $0$, $1$ or $2$. In order to show that $\Hom(E_2,E_0[k])=0$, notice that, by Proposition III.6.7 in \cite{Har}, these are isomorphic to $\Ext^k(\mathcal{O}_{X},2L')$, that is, to the cohomologies $H^k({X};2L')$. On the other hand, these all vanish by the first part of the proof. Next, the vanishings of the $\Hom(E_2,E_1[k])$'s and the $\Hom(E_1,E_0[k])$'s are equivalent to those of the $\Ext^k(\mathcal{O}_X,L')\simeq H^k(X;L')$ above. We are left with showing that $\Hom(E_i,E_i[k])$, for $i=0,1,2$, has rank one for $k=0$ and vanishes for higher values of $k$. In all of the above instances, we have $\Hom(E_i,E_i[k])\simeq \Hom(E_0,E_0[k])\simeq H^k(X;\mathcal{O}_X)$, which has rank one for $k=0$ and vanishes for higher values of $k$, since any fake projective plane has the same Hodge diamond of a projective plane. Conversely, notice that the implications in the above argument can be reversed. This concludes our proof.
\end{pf}

The next proposition gives a vanishing result for certain $\Z/7\Z$-equivariant line bundles on Keum's fake projective planes. This vanishing plays a crucial role in the proof of Theorem \ref{second}. 

\begin{prop}\label{vanishing}
Let $X$ be a fake projective plane with $Aut(X)=G_{21}$. Let $\Z/7\Z\lhd G_{21}$ be the unique $7$-Sylow. If the line bundle $2\mathcal{O}_{X}(1)+T$ is $\Z/7\Z$-equivariant for some torsion line bundle $T$, we have $H^{0}(X; 2\mathcal{O}_{X}(1)+T)=0$.
\end{prop}
\begin{pf}
By contradiction, let us assume $H^{0}(X; 2\mathcal{O}_{X}(1)+T)\neq 0$. Since the line bundle $2\mathcal{O}_{X}(1)+T$ is $\Z/7\Z$-equivariant, we have that $\Z/7\Z$ acts on the projectivization of $H^{0}(X; 2\mathcal{O}_{X}(1)+T)$. By Lemma \ref{L2}, this projective space can either be $\pp^{0}$ or $\pp^{1}$. Since any finite group action on $\pp^{1}$ has necessarily a fixed point, we have the existence of a curve in the linear system $|2\mathcal{O}_{X}(1)+T|$, say $C$, which is $\Z/7\Z$-invariant. Next, we claim that such a curve cannot be pointwise fixed by the $\Z/7\Z$-action. This nontrivial fact follows from the ground breaking work of Keum. More precisely, Proposition 2.4 and Theorem 1.1. in \cite{Keum} imply that the dimension of the fixed point set of a  $\Z/7\Z$-action on a fake projective plane is zero. In particular, there cannot be any curve in $X$ pointwise fixed by the $\Z/7\Z$-action. Moreover, since on Keum's fake projective plane there are no curves linearly equivalent to $L_{1}$, the curve $C$ fixed by the $\Z/7\Z$-action is necessarily reduced and irreducible.

Thus, we have a $\Z/7\Z$-invariant curve $C$ such that $p_{a}(C)=6$. If $C$ is smooth we have $p_{a}(C)=g(C)=6$. Now a nontrivial $\Z/7\Z$-action on a genus $6$ curve has necessarily $4$ fixed points, see Proposition V.2.14 page 285 in \cite{Kra}. This contradicts the fact that $X/(\Z/7\Z)$ has exactly $3$ singular points, see again Theorem 1.1 in \cite{Keum}. We therefore conclude that $C$ cannot be smooth and let us denote by $\overline{C}$ its normalization. Now, let us observe that the nontrivial $\Z/7\Z$-action on $C$ lifts to the normalization $\overline{C}$. Thus, $\overline{C}$ is a smooth curve of genus strictly less than $6$ with a nontrivial automorphism of order $7$. By Proposition V.2.14. page 285 in \cite{Kra}, we conclude that the genus of $\overline{C}$ must necessarily be equal to $3$. The theory of the Toledo invariant applied to ball quotient surfaces, see for example Proposition 2.1 in \cite{ultimo}, now tells us that
\begin{align}\label{Domingo}
K_{X}\cdot C\leq 3(g(\overline{C})-1)
\end{align}
with equality if and only if $C$ is a totally geodesic immersed curve. Thus, since we established $g(\overline{C})=3$ and 
\[
K_{X}\cdot C=3L_{1}\cdot 2L_{1}=6
\]
we can use (\ref{Domingo}) to deduce that $C$ is necessarily a totally geodesic immersed curve.  On the other hand it is well known that there are no immersed totally geodesic curves in any of the fake projective planes. In fact, the main theorem in \cite{Prasad} combined with \cite{Steger} tells us that any fake projective plane is an arithmetic ball quotient of the second type. Now, ball quotients of the second type cannot contain any immersed totally geodesic curve, see for example page 901 in \cite{Moller-Toledo}. The proof is therefore complete.
\end{pf}

We are now ready to produce a nonstandard exceptional collection on at least one of Keum's fake projective planes.

\begin{thm}\label{eccezionale}
Let $X$ be the fake projective plane with $Aut(X)=G_{21}$ and $H_{1}(X; \Z)=(\Z/2\Z)^{4}$. There exist a nontrivial torsion line bundle $T\in Pic(X)$, such that the collection
\[
\mathcal{O}_{X}, -\mathcal{O}_{X}(1)-T, -2\mathcal{O}_{X}(1)
\] 
is exceptional.
\end{thm}
\begin{pf}
Recall that the standard exceptional collection 
\[
\mathcal{O}_{X}, -\mathcal{O}_{X}(1), -2\mathcal{O}_{X}(1)
\]
was constructed by Galkin \emph{et al}. in \cite{Katzarkov}. In order to construct the nonstandard collection, by Proposition \ref{P:cubicroot} it is enough to show that there exists a numerical cubic root of $K_{X}$ \emph{not} linearly equivalent to $\mathcal{O}_{X}(1)$, say $L'$, for which the cohomologies $H^{0}(X; 2L')$, $H^{2}(X; L')$ and $H^{2}(X; 2L')$ all vanish. First, recall that we have $2L'\cong 2\mathcal{O}_{X}(1)$, thus the vanishing of $H^{0}(X; 2L')$ follows directly from Theorem 1.3 in \emph{loc. cit.}. Next, observe the following. Given any ample generator of Pic(X), say $L'$, we can write $L'\cong \mathcal{O}_{X}(1)+T_{i}$, where $T_{i}$ is a torsion line bundle. Thus, if $X$ is any of Keum's fake projective planes, by duality we have
\begin{align}\notag
H^{2}(X; 2L')& = H^{0}(X; K_{X}-2L')= H^{0}(X; 3\mathcal{O}_{X}(1)-2\mathcal{O}_{X}(1)-2T_{i})\\ \notag
&=H^{0}(X; \mathcal{O}_{X}(1))=0,
\end{align}
since $2T_{i}\cong \mathcal{O}_{X}$. Next, let us observe that
\begin{align}\notag
H^{2}(X; L')&=H^{0}(X; K_{X}-L')=H^{0}(X; 3\mathcal{O}_{X}(1)-\mathcal{O}_{X}(1)-T_{i}) \\ \notag
&=H^{0}(X; 2\mathcal{O}_{X}(1)+T_{i}),
\end{align} 
since any torsion line bundle $T_{i}\in Pic(X)$ has order two.
In conclusion, if we can find a line bundle $2\mathcal{O}_{X}(1)+T_{i}$ with $T_{i}\neq 0$ and such that $H^{0}(X; 2\mathcal{O}_{X}(1)+T_{i})=0$, we have that the collection
\begin{align}\label{exceptional}
\mathcal{O}_{X}, -\mathcal{O}_{X}(1)-T_{i}, -2\mathcal{O}_{X}(1)
\end{align}
is exceptional. Now, let $X$ be the unique fake projective plane with $H_{1}(X; \Z)=(\Z/2\Z)^{4}$, so that there are exactly fifteen nontrivial torsion line bundles on $X$. Thus, there is at least one nontrivial torsion line bundle, say $T$, which is $\Z/7\Z$-equivariant. Here the $\Z/7\Z$-action comes from the unique $7$-Sylow in $Aut(X)=G_{21}$. Now, if $T$ is a nontrivial $\Z/7\Z$-equivariant torsion line bundle on $X$, we then claim that $H^{0}(X; 2\mathcal{O}_{X}(1)+T)=0$, which is enough to show that the collection given in \eqref{exceptional} is indeed exceptional. In order to show that $H^{0}(X; 2\mathcal{O}_{X}(1)+T)=0$, let us observe that the line bundle $2\mathcal{O}_{X}(1)+T$ is $\Z/7\Z$-equivariant. In fact, $2\mathcal{O}_{X}(1)$ is $G_{21}$-equivariant and the torsion part $T$ is by construction $\Z/7\Z$-equivariant.  It then suffices to apply Proposition \ref{vanishing} to conclude the proof. In conclusion, there exists a nontrivial torsion line bundle $T$ such that the collection
\begin{equation}\label{e:exotica}
\mathcal{O}_{X}, -\mathcal{O}_{X}(1)-T, -2\mathcal{O}_{X}(1)
\end{equation}
is exceptional.
\end{pf}

\begin{rem}\label{R:exotica}
In the terminology of \cite{Galkin}, a $2$-\emph{minifold} is a two dimensional smooth complex projective variety $X$ admitting a full exceptional collection of length $3$ in its derived category $\mathcal{D}^b(X)$. Thanks to the classification of minifolds appearing in Theorem 1.1 of \emph{loc. cit.}, the only $2$-{minifold} is $\mathbb{P}^2$, hence the exceptional collection from \eqref{e:exotica} cannot be full. In particular, it makes sense to consider its right orthogonal in $\mathcal{D}^b(X)$.
\end{rem}

The following corollary to our main theorem tells us that we have obtained yet another semi-orthogonal decomposition of the bounded derived category of coherent sheaves $\mathcal{D}^b(X)$.

\begin{cor}\label{orthogonal}
Let $X$ be as in Theorem \ref{eccezionale} and let $\mathcal{A}$ be a right orthogonal of the exceptional collection from \eqref{e:exotica}. Then there is a semi-orthogonal decomposition
$$
\mathcal{D}^b(X)=\langle \mathcal{O}_{X}, -\mathcal{O}_{X}(1)-T, -2\mathcal{O}_{X}(1), \mathcal{A}\rangle.
$$
\end{cor}

\begin{pf}
Again, this is done along the lines of the proof of Corollary 1.2 in \cite{Katzarkov}. Let $\mathcal{B}$ denote the triangulated envelope $\langle \mathcal{O}_{X}, -\mathcal{O}_{X}(1)-T, -2\mathcal{O}_{X}(1)\rangle$ in $\mathcal{D}^b(X)$. Since we proved above that \eqref{e:exotica} is an exceptional collection, the category $\mathcal{B}$ is (left- and right-) admissible by Theorem 3.2a in \cite{Bondal}. Now, since we defined $\mathcal{A}$ as $\mathcal{B}^{\perp}$, from Lemma 3.1 in \emph{loc. cit.} it follows that $\mathcal{B}$ and $\mathcal{A}$ generate $\mathcal{D}^b(X)$.
\end{pf}

\begin{rem}\label{R:phantom}
One might ask whether our $X$ is such that its derived category $\mathcal{D}^b(X)$ admits an $H$-phantom. This question has an affirmative answer. In fact, in our case the canonical class $K_X$ is divisible by $3=dim(X)+1$, hence the existence of an $H$-phantom is implied by Corollary 1.2 in \cite{Katzarkov}. More in detail, the subcategory $\mathcal{A}$ in the statement is indeed an $H$-phantom. In order to see this, one can follow Section 4 in \cite{Keumvan} and recall that the Hochschild \emph{homology} of $X$ is isomorphic to its Hodge cohomology. Since $X$ is a fake projective plane, its Hodge cohomology must have dimension equal to $3$. Now, since the sum of the Hodge numbers of our exceptional collection is already $3$ and Hochschild homology is additive on semiorthogonal decompositions (see Section 7.2 in \cite{Kuz}), we must have that $HH_{\bullet}(\mathcal{A})=0$.
\end{rem}


\end{document}